\documentclass{amsart}
\usepackage[cp1251]{inputenc}
\usepackage[english,russian]{babel}
\usepackage{amsmath,graphicx}
\usepackage{amssymb}
\usepackage{amsfonts}

\begin{document}
\thispagestyle{empty}
\title[Second order ODE's \dots ]
{Second order ODE's cubic in the first order derivative 
with 2-dimensional symmetry algebra}

\author{ Vera V. Kartak~$^\dag$}

\address{$^\dag$~450000, Ufa State Aviation Technical University, Ufa, K.Marx str., 12, Russia; 450076, Bashkir State University, Ufa, Z.Validi str., 32, Russia}
\email{kvera@mail.ru} 


\maketitle 
{
\small
\begin{quote}
\noindent{\bf Abstract. } 
We describe the second order ODE's cubic in the first order derivative with 2-dimensional symmetry algebra. We show that there exist only eight different types of them. We also construct the easily verifiable Equivalence Criterion for every type of equations analogous to Linearization Criterion of S. Lie.
\end{quote}



\section{Introduction}

The group classification
of the second order ODEs of the following form
\begin{equation}\label{li}
y''=f(x,\,y,\,y')
\end{equation}
was provided in the works of S.Lie, see \cite{Lie}, A.Tresse, see \cite{Tresse2}. 
More details see in modern works of N.Ibragimov \cite{Ibragimov}, B.Kruglikov \cite{Kruglikov}.  
It turns out that the dimension of symmetry algebra is lacunary, it can be only 8, 3, 2, 1 or 0. 

For the equations with 8-dimensional symmetry algebra it is well-known the Linearization criterion, see
\cite{Cartan}, \cite{Ibr1}, \cite{Grissom}, \cite{Gonz}, \cite{Qadir}, \cite{Kamran}, \cite{Tresse2}, \cite{Yum}, \cite{Yum1}, \cite{Rom}, \cite{Sharipov2}.

 {\bf Linearization Criterion.} {\it These propositions are equivalent:
 \begin{enumerate}
 \item Equation (\ref{li}) has a 8-dimensional symmetry algebra;
 \item equation  (\ref{li}) is linearizable, by the generic point transformation 
  \begin{equation}\label{zam}
 \tilde x=\tilde x(x,\, y),\qquad \tilde y=\tilde y(x,\, y)
 \end{equation}  
 it takes the form $\tilde y''=0$;
 \item equation (\ref{li}) has form
 \begin{equation}\label{eq}
y''=P(x,y)+3\,Q(x,y)y'+3\,R(x,y)y^{\prime 2}+S(x,y)y^{\prime 3}
\end{equation}
  and the following conditions are true
 \begin{equation}\label{alpha}
\aligned&\aligned A=P_{ 0.2}&-2Q_{ 1.1}+R_{ 2.0}+ 2PS_{ 1.0}+SP_{
1.0}-\\
&-3PR_{ 0.1}-3RP_{ 0.1} -3QR_{ 1.0} +6QQ_{ 0.1}=0,\endaligned \\
\vspace{1ex}&\aligned B=S_{ 2.0}&-2R_{ 1.1}+Q_{ 0.2}- 2SP_{
0.1}-PS_{ 0.1}+\\
&+3SQ_{ 1.0}+3QS_{ 1.0}+ 3RQ_{ 0.1}-6RR_{ 1.0}=0.
\endaligned\endaligned
\end{equation}
(Here and below $F_{i.j}=\partial ^{i+j}F/\partial x^i\partial y^j$).
\end{enumerate}}

As we can see for the arbitrary equation (\ref{eq}) it is very easy to check the Lie's Linearization Criterion. However usually it is a non-trivial problem to find the corresponding change of variables (\ref{zam}) that reduces our equation into the form $\tilde y''=0$, see \cite{Kruglikov}.

{\bf Example.} The equations No.  6.113, 6.134, 6.169 in the Handbook by E.Kamke \cite{Kamke} satisfy the conditions of Linearization Criterion. 
$$
\aligned
6.113\quad & yy''-y'^2-y^2\ln y=0,\\
6.134\quad & (y-x)y''-2y'(y'+1)=0,\\
6.169\quad & xyy''+xy'^2-yy'=0.
\endaligned
$$

Equations (\ref{li}) with 3-dimensional symmetry algebra also investigated in detail, see \cite{Tresse2},
\cite{Rom}, \cite{Ibragimov},  \cite{Sharipov2}, \cite{IbrMel1}, \cite{IbrMel2}.
By the point transformations (\ref{zam}) they could be reduced into the one of the canonical (or normal) form
$$
\begin{aligned}
y''={y'}^a,&\qquad\qquad y''=\frac{(cy'+\sqrt{1-y'^2})(1-y'^2)}{x},\\
y''=e^{y'},&\qquad\qquad y''=\pm(xy'-y)^3.
\end{aligned}
$$
(In the last formulas we should not write  tildas on $x$, $y$, $y'$ and $y''$.)

Substantial progress is in the equations (\ref{li}) with 2-dimensional symmetry algebra.
According to the S.Lie \cite{Lie}, see also \cite{Ibragimov}, \cite{Kruglikov},
 all equations (\ref{li}) addmitting a 2-dimensional subalgebra could be reduced
by the transformations (\ref{zam}) into the one of the following forms, here $f(x)$ and $f(y')$  is the certain arbitrary functions
\begin{equation}\label{types}
y''=f(y'),\qquad y''=f(x),\qquad y''=\frac 1x f(y'),\qquad y''=f(x)y'.
\end{equation}
Let us note that the equations $y''=f(x)$ and $y''=f(x)y'$ from the list (\ref{types}) satisfy the Linearization Criterion.

For the equations (\ref{li}) in the particular form $y''=f(x,y)$ with 2-dimensional subalgebra was made the more precise classification, see S.Lie  \cite{Lie} and  L.Ovsyannikov  \cite{Ovsyan}. 
The complete list of the canonical forms of these equations is the following
\begin{equation}\label{ovs}
y''=e^y,\quad\qquad y''=y^k,\;k\ne -3,\quad\qquad y''=\frac 1{y^3}.
\end{equation}
What is more, the equation  $y''= 1/{y^3}$ has a 3-dimensional symmetries algebra.

In the paper \cite{Ovsyan} L.Ovsyannikov wrote:  "But formed in this papers analysis does not work
the complete solution of the equivalence problem consisting in 
establishing {\it criteria equivalence for the a-priori given equations} $y''=f(x,y)$ and $\tilde y''=\tilde f(\tilde x,\tilde y)$. Such a criterion may be obtained by only on the basis of the theory of differential invariants."

The Invariant Theory of the second order ODEs has been extensively studies. See works \cite{Lie1}, \cite{Lie}, \cite{Liouville}, \cite{Tresse1}, \cite{Tresse2}, \cite{Cartan}, \cite{Kamran}, \cite{Thomsen},
\cite{Sokolov}, \cite{Milson}, \cite{Kruglikov}, \cite{Sharipov1}, \cite{Sharipov2}, \cite{BordagBandle}, etc.

In the paper \cite{Sokolov} was set up the Problem: {\it "Find necessary and sufficient conditions for a given equation of the form (\ref{li}) to be point-equivalent to one of the model equations, admitting a two-dimensional Lie group of point symmetries". }
 
In the work \cite{Kamran} was made the complete classification of second order ODEs addmitting Lie groups of {\it fibre-preserving point symmetries}. It means that were considered only restricted transformations:
$$
  \tilde x=\tilde x(x),\qquad \tilde y=\tilde y(x,\, y).
$$
We also note the following interesting results. In the work \cite{Mahomed} were founded the five representations of equivalence classes of equations (\ref{li}) (as the model equations) and in the work \cite{Nest} were investigated the finite-dimensional Lie groups acting on the real plane.

At the present paper this Problem completely solved not for the general form of equations (\ref{li}) but only for the equations of the form (\ref{eq}). Were founded eight different types of equations (\ref{eq}) with 2-dimensional subalgebra that are the model equations (including 8- and 3-dimensional cases). 
For every model equation the necessary and sufficient conditions of equivalence (the Equivalence Criterion) are constructed. They are similar like to the Linearization Criterion of S.Lie.

\newpage
\section{The main results}
There are exist only eight different types of equations (\ref{eq}) with 2-dimensional subalgebra.
A full list of them is presented into the Table. Equations from the different lines of table are not point-equivalent. Equations that are from the one line of Table point equivalent if and only if all invariants are coincide ($k$ in the lines 4, 5, 7 and $n$, $a$, $b$ in the line 8).

\begin{table}[htbp]
\centering
\begin{tabular}{|c|l|l|l|}
\hline \hline   
$Dim (Z)$ & Equivalence & Model Equation & Algebra \\
 & Criterion &  & $X=\xi\frac{\partial}{\partial x}+\eta\frac{\partial}{\partial y}$ \\
\hline 
 & &  & $\xi=0,\,\eta=1;$ \\
 8 &   Linearization            &  $y''=0$       & $\xi=1,\, \eta=0;$\\
  &   criterion            &         & $\xi=0,\, \eta=x;$\\
  &              &         & $\xi=0,\, \eta=y;$\\
  &              &         & $\xi=x,\, \eta=0;$\\
  &              &         & $\xi=y,\, \eta=0;$\\
  &              &         & $\xi=xy,\, \eta=y^2;$\\
  &              &         & $\xi=x^2,\, \eta=xy$\\
\hline 
  &  &  & $\xi=1,\,\eta=0;$ \\
  3 & Theorem 1 & $y''=\frac 1{y^3}$ & $\xi=2x,\,\eta=y;$ \\
   &   &  & $\xi=x^2,\,\eta=xy$ \\ 
\hline & & & \\
  &  &  & $\xi=1,\,\eta=0;$ \\ 
 2  & Theorem 2 & $y''=e^y$ & $\xi= x,\,\eta=-2$ \\
  &  &  &  \\
\hline    
  & Theorem 3 &  & $\xi=1,\,\eta=0;$ \\
 2  & $k$ -  invariant & $y''=\frac{y^{k+2}}{(k+1)(k+2)}$ & $\xi={x(k+1)},\,\eta=-2y$ \\
   & $k\ne 0,\, -1,\, -2,\, -5$ &  &  \\
  \hline   & & & \\
 & Theorem 4 &  & $\xi=1,\,\eta=0;$ \\
 2  & $k$ -  invariant & $y''=\frac{y'^2}{2y}+\sqrt{6yk}y'+$  & $\xi= x,\,\eta=-2y$ \\
   & $k\ne 0$ & $\phantom{....}+\frac{2(1-2k)y^2}3$ &  \\
  \hline   
& & & \\
2 &  Theorem 5 & $y''=\frac{y^2}2$ & $\xi=1,\,\eta=0;$ \\
  & &  & $\xi= x,\,\eta=-2y$ \\ 
\hline  & & & \\  
 2 & Theorem 6 & $y''=y'^3-ky y'^2+\frac{k^2y^2y'}3+$ & $\xi=1,\,\eta=0;$ \\
  & $k$ -  invariant   & $\phantom{....}+\frac 1{k}+\frac{k^2y}9-\frac{k^3y^3}{27} $ & $\xi=0,\,\eta=e^{kx/3}$ \\ 
 & $k\ne 0$ & &\\
\hline    
& & & \\
 & Theorem 7 & $y''=\frac {4y'^3}{ny^3}+\frac{((6a-1)y+6)y'^2}{y^{2}}+$ & $\xi=1,\,\eta=0;$ \\
 2 & $ c,\, d,\, k,\, m$ from (\ref{coeff}) & $\phantom{....}+\frac{3(by^2+cy+n)y'}{y}+ $ & $\xi=e^{nx/2},$ \\
     & $n,\,a,\, b$ -  invariants & $\phantom{....}+\frac 12(dy^3+ky^2+my+n^2)$ & $\eta=-\frac{ny}2e^{nx/2}$ \\
   &$n\ne 0$ & & \\
\hline\hline    
\end{tabular}
\end{table}

The constants  $c$, $d$, $k$, $m$ of the model equation from the line 8 of the Table
are expressed into the invariants
 $n$, $a$ and $b$ by the following formulas:
 \begin{equation}\label{coeff}
 \aligned
  c& =\frac{n(12a-5)}{6}, \quad d =\frac{bn(6a-1)}{6},\quad
k =3bn,\quad m=\frac{3n^2(2a-1)}2,\\
&\text{where}\qquad 18b(b-a^2n)+15nab-2bn-6=0.
 \endaligned
\end{equation}

\section{The method of solution}

The method of the constructing the invariants of the equations (\ref{eq})
based on the works \cite{Sharipov1}-\cite{Sharipov3}, \cite{Kartak3}. It is a geometrical approach
that allows us to get the explicit formulas for the invariants.

With the equations (\ref{eq}) are associated the pseudovectorial fields:
${\boldsymbol \alpha}=(B,\, -A)$ with the components from (\ref{alpha}) and ${\boldsymbol \beta}=(G,\, H)$.
The prefix 'pseudo'  means that under the change of variables the fields ${\boldsymbol \alpha}$ and ${\boldsymbol \beta}$ transform by the rule of vectorial fields with the factor
$$ \aligned G&=-BB_{ 1.0}-3AB_{ 0.1}+4BA_{ 0.1}+
3SA^2-6RBA+3QB^2,\\
H&=-AA_{ 0.1}-3BA_{ 1.0}+4AB_{ 1.0}-
3PB^2+6QAB-3RA^2.
\endaligned
$$

The function $F$ is defined via the scalar product of the  fields ${\boldsymbol \alpha}$ and ${\boldsymbol \beta}$ using by skewsymmetric Gramian matrix
\begin{equation}\label{F}
3F^5=({\boldsymbol \alpha},\, {\boldsymbol \beta})=d_{ij}\alpha^i\beta^j=AG+BH,\qquad d_{ij}=\left(\begin{array}{rr}  0 & 1\\ -1 & 0
\end{array}\right).
\end{equation}
Details see in works \cite{Sharipov1}-\cite{Sharipov3}. 
There are only three possible cases:
\begin{enumerate}
\item ${\boldsymbol \alpha}=0$, it means $A=0$ и $B=0$ from (\ref{alpha}), {\it Maximal Degeneration Case}. This case corresponds to Lie's Linearization Criterion;
\item ${\boldsymbol \alpha}\ne 0$ and the fields ${\boldsymbol \alpha}$ and ${\boldsymbol \beta}$ are collinear, $({\boldsymbol \alpha},\, {\boldsymbol \beta})=0$, it means that $A\ne 0$ or $B\ne 0$ from (\ref{alpha}), but $F=0$ from (\ref{F}), {\it Intermediate Degeneration Case};
\item the fields ${\boldsymbol \alpha}$ and ${\boldsymbol \beta}$ are non-collinear, $({\boldsymbol \alpha},\, {\boldsymbol \beta})\ne 0$, it means that $F\ne 0$ from (\ref{F}), {\it General Case}.
\end{enumerate}
Intermediate Degeneration Case splits into seven Cases. The First and the Seventh Cases split into Subcases 1.1, 1.2, 1.3, 1,4, 7.1, 7.2. 

The following diagram illustrates the Intermediate Degeneration  Сases.
\begin{figure}[htbr]
\center{\includegraphics[width=1\linewidth]{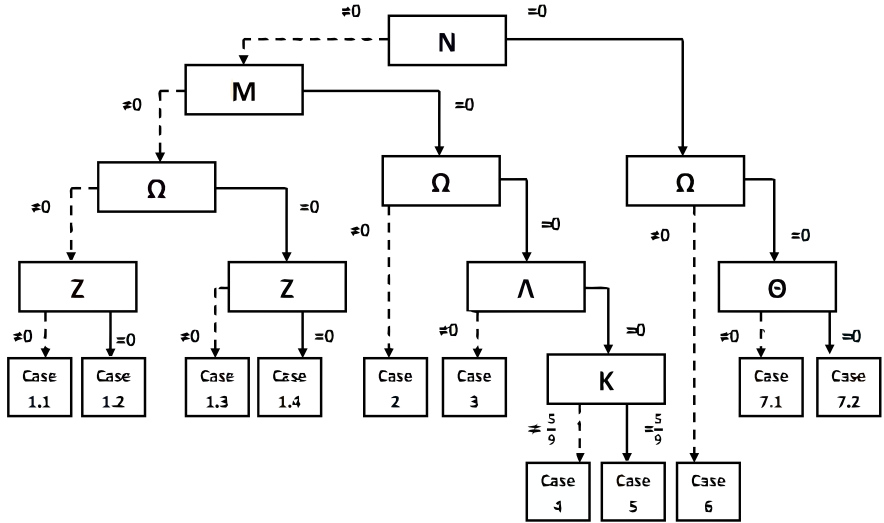}}
\caption{Tree of  Intermediate Degeneration Cases.}
\label{tree}
\end{figure}

Here $M$, $N$, $\Omega$, $Z$, $\Lambda$, $K$, $\Theta$ -- pseudoinvariants --  the functions depending on $(x,\,y)$ such that under the transformations (\ref{zam}) they transforms with the weight: $\tilde J=\mu(x,y)J$. 
If for the certain equation (\ref{eq}) the pseudoinvariant is vanishes, then it becomes to be an invariant. All given (pseudo)invariants are calculated via the coeffitions of the equation (\ref{eq}) on the direct formulas. These formulas will be  given below.

 In the paper \cite{Sharipov2} was proved the following Proposition.
 
 {\bf Proposition. }{\it Equations (\ref{eq}) with the 8-dimensional  symmetry algebra 
 are only in the Case of Maximal Degeneration; with the 3-dimensional  symmetry algebra
 are only in the Fifth Case of Intermediate Degeneration; with the 2-dimensional  symmetry algebra only in the cases:
\begin{enumerate}
\item General Case when all invariants are the constants;
\item First Case of Intermediate Degeneration when all invariants are the constants;
\item Seventh Case  of Intermediate Degeneration when the pseudionvariant $L=0$ (\ref{L}).
\end{enumerate}
 }

 Let us describe in detail all this cases.
 
 \section{Equations (\ref{eq}) with the 3-dimensional symmetry algebra}
 In this case we repeat the results from the works \cite{Rom}, \cite{Sharipov2}, \cite{Sharipov3}, \cite{Ovsyan}.
 
As the model equations proposed the following equations:
$$
\text{\cite{Rom}:}\quad y''=\pm (xy'-y)^3,\quad\text{\cite{Sharipov2}, \cite{Sharipov3}:}\quad y''=-\frac{5}{4}y'+\frac 43 x^2y^{\prime 3},\quad\text{\cite{Ovsyan}:}\quad y''=\frac{1}{y^3}.
$$

{\bf Theorem 1. }
 {\it These propositions are equivalent:
 \begin{enumerate}
 \item Equation (\ref{eq}) has the 3-dimensional  point symmetries algebra;
 \item Equation (\ref{eq})  reduced to the form $y''= y^{-3}$ by the point transformations (\ref{zam});
 \item Equation (\ref{eq}) is from the Fifth Case of the Intermediate Degeneration: $F=0$ from (\ref{F}), $A\ne 0$ or $B\ne 0$ from (\ref{alpha}), $N\ne 0$ from (\ref{N}), $M=0$ from (\ref{M1}), (\ref{M2}), $\Omega=0$ from (\ref{Omega1}), (\ref{Omega2}), $\Lambda=0$ from (\ref{lambda}), $K=- 5/9$  from (\ref{K}).
 \end{enumerate}}

{\bf Example.} The equations No.  6.81, 6.138 in the Handbook by E.Kamke \cite{Kamke} satisfy the conditions of Theorem 1
$$
\aligned
6.81\quad & 2xy''+y'^3+y'=0,\\
6.138\quad & 2yy''-y'^2+a=0,\quad a=const.
\endaligned
$$

In the cases $A\ne 0$ and $B\ne 0$ the pseudoinvariant $N$ is given by the formulas
\begin{equation}\label{N}
 N =-\frac H{3A}, \qquad\qquad N =\frac G{3B}.
\end{equation}
If at the same time $A\ne 0$ and $B\ne 0$ we can use any of the proposed formulas.

The pseudoinvariant $M$  in the case $A\ne 0$ reads as
\begin{equation}\label{M1}
\aligned
M=&-\frac {12BN(BP+A_{ 1.0})}{5A}+BN_{ 1.0}+\frac {24}5BNQ+\\
&+\frac 65NB_{ 1.0}+\frac 65NA_{ 0.1}-AN_{ 0.1}- \frac {12}5ANR
\endaligned
\end{equation}
and in the case $B\ne 0$ it is given by the formula
\begin{equation}\label{M2}
\aligned
M=&-\frac {12AN(AS-B_{ 0.1})}{5B}-AN_{ 0.1}+\frac {24}5ANR-\\
 &-\frac
65NA_{ 0.1}-\frac 65NB_{ 1.0}+BN_{ 1.0}- \frac {12}5 BNQ.
\endaligned
\end{equation}

As $A\ne 0$, the explicit formula for the pseudoinvariant  $\Omega$ reads as
\begin{equation}\label{Omega1}
\aligned 
\Omega &=\frac {2BA_{ 1.0}(BP+ A_{ 1.0})}{A^3}- \frac
{(2B_{ 1.0}+3BQ)A_{ 1.0}}{A^2}+\frac {(A_{ 0.1}-2B_{1.0})BP}{A^2}+ 
\frac {B_{ 2.0}}A-\\
&-\frac {BA_{ 2.0}+B^2 P_{ 1.0}}{A^2}+\frac {3B_{ 1.0}Q+3BQ_{ 1.0}- B_{ 0.1}P-BP_{
0.1}}{A}+Q_{ 0.1}- 2R_{ 1.0}.
\endaligned
\end{equation}

And in the case  $B\ne 0$ the similar formula is
\begin{equation}\label{Omega2}
\aligned 
\Omega &=\frac {2AB_{ 0.1}(AS- B_{ 0.1})}{B^3}- \frac
{(2A_{ 0.1}-3AR)B_{ 0.1}}{B^2}+\frac {(B_{ 1.0}-2A_{0.1})AS}{B^2}- 
\frac {A_{ 0.2}}B+ \\
&+\frac {AB_{ 0.2}-A^2 S_{ 0.1}}{B^2}+\frac {3A_{ 0.1}R+3AR_{ 0.1}- A_{ 1.0}S-AS_{
1.0}}{B}+R_{ 1.0}- 2Q_{ 0.1}.
\endaligned
\end{equation}

The pseudoinvariant $\Lambda$ in the cases $A\ne 0$ and $B\ne 0$ respectively reads as
\begin{equation}\label{lambda}
\aligned
\Lambda=&\frac{6N(BP+B_{1.0})}{5A^2}-\frac{N_{1.0}}A-\frac{6NQ}{5A}-2\Omega,\\
\Lambda=&-\frac{6N(AS-B_{0.1})}{5B^2}-\frac{N_{0.1}}B+\frac{6NR}{5B}-2\Omega;
\endaligned
\end{equation}
the pseudocovectorial field $\omega$ in the case $A\ne 0$ is given by the formula
\begin{equation}\label{omega1}
\begin{aligned}
\omega_1&=\frac{12PR}{5A}-\frac{54}{25}\frac{Q^2}{A}-\frac{P_{0.1}}{A}+\frac{6Q_{1.0}}{5A}-
\frac{PA_{0.1}+BP_{1.0}+A_{2.0}}{5A^2}-\\
&-\frac{2B_{1.0}P}{5A^2}+\frac{3QA_{1.0}-12PBQ}{25A^2}+
\frac{6B^2P^2+12A_{1.0}BP+6A_{1.0}^2}{25A^3},\\
\omega_2&=\frac{6\Lambda+3\Omega}{5A}+\frac{-5BP_{0.1}+6BQ_{1.0}+12RBP}{5A^2}-
\frac{54}{25}\frac{BQ^2}{A^2}-\frac{12B^2PQ-3BQA_{1.0}}{25A^3}-\\
&-\frac{2BB_{1.0}P+BA_{0.1}P+B^2P_{1.0}+BA_{2.0}}{5A^3}+\frac{6BA_{1.0}^2+6B^3P^2+12B^2A_{1.0}P}{25A^4},
\end{aligned}
\end{equation}
in the case $B\ne 0$ the similar formula is
\begin{equation}\label{omega2}
\begin{aligned}
\omega_1&=-\frac{6\Lambda+3\Omega}{5B}+\frac{5AS_{1.0}-6AR_{0.1}+12QAS}{5B^2}-
\frac{54}{25}\frac{AR^2}{B^2}-\frac{12A^2SR-3ARB_{0.1}}{25B^3}+\\
&+\frac{2AA_{0.1}S+AB_{1.0}S+A^2S_{0.1}-AB_{0.2}}{5B^3}+\frac{6AB_{0.1}^2+6A^3S^2-12A^2B_{0.1}S}{25B^4},\\
\omega_2&=\frac{12SQ}{5B}-\frac{54}{25}\frac{R^2}{B}+\frac{S_{1.0}}{B}-\frac{6R_{0.1}}{5B}+
\frac{SB_{1.0}+AS_{0.1}-B_{0.2}}{5B^2}+\\
&+\frac{2A_{0.1}S}{5B^2}-\frac{3RB_{0.1}+12SAR}{25B^2}+
\frac{6A^2S^2-12B_{0.1}AS+6B_{0.1}^2}{25B^3};
\end{aligned}
\end{equation}
the pseudoinvariant $K$  in the cases $A\ne 0$ and $B\ne 0$ respectively reads as
\begin{equation}\label{K}
\aligned
K=&\frac{\Lambda_{1.0}+\Lambda\varphi_1}{A}+\frac{\Omega_{1.0}+\Omega\varphi_1}{3A}+
\frac{N\omega_1}{A},\qquad A\ne 0,\\
K=&\frac{\Lambda_{0.1}+\Lambda\varphi_2}{B}+\frac{\Omega_{0.1}+\Omega\varphi_2}{3B}+
\frac{N\omega_2}{B},\qquad B\ne 0,
\endaligned
\end{equation}
where $\varphi_i$ in the cases $A\ne 0$ and $B\ne 0$ respectively are given by the formulas
\begin{equation}\label{phi}
\aligned
\varphi_1=&-3\frac {BP+A_{ 1.0}}{5A}+\frac 35Q, \;\;
\varphi_2=3B\frac {BP+A_{ 1.0}}{5A^2} -3\frac {B_{ 1.0}+A_{
0.1}+3BQ}{5A}+\frac 65R,\\
 \varphi_1=&-3A\frac {AS-B_{ 0.1}}{5B^2}
-3\frac {A_{ 0.1}+B_{ 1.0}-3AR}{5B}-\frac 65Q,\quad
\varphi_2=3\frac {AS-B_{ 0.1}}{5B}-\frac 35R.
\endaligned
\end{equation}

 \section{Equations (\ref{eq}) with the 2-dimensional  symmetry algebra}

\subsection{First Case of the Intermediate Degeneration}

In this case the basic invariants of the equation (\ref{eq}) are calculated by the formulas
\begin{equation}\label{inv1}
I_1=\frac{M}{N^2},\quad I_2=\frac{\Omega^2}{N},\quad I_3=\frac{\Gamma}{M},
\end{equation}
where $N$ from (\ref{N}), $M$ from (\ref{M1}), (\ref{M2}), $\Omega$ from (\ref{Omega1}), (\ref{Omega2}), $\Gamma$ from (\ref{Gamma}).
\begin{equation}\label{Gamma}
\aligned
\Gamma=&\frac {\gamma^1\gamma^2(\gamma^1_{
1.0}- \gamma^2_{ 0.1})}{M}+ \frac {(\gamma^2)^2\gamma^1_{ 0.1}-
(\gamma^1)^2\gamma^2_{ 1.0}}M+\\
&+\frac
{P(\gamma^1)^3+3Q(\gamma^1)^2\gamma^2+3R\gamma^1(\gamma^2)^2+
S(\gamma^2)^3}M,
\endaligned
\end{equation}
here the functions $P$, $Q$, $R$, $S$ are the coefficients of the equation  (\ref{eq}), pseudovectorial field $\gamma=(\gamma^1,\,\gamma^2)$ in the case $A\ne 0$ reads as
\begin{equation}\label{gamma1}
\aligned
 \gamma^1=&-\frac {6BN(BP+A_{ 1.0})}{5A^2}+
\frac {18NBQ}{5A}+
\frac {6N(B_{ 1.0}+A_{ 0.1})}{5A} -N_{ 0.1}-\frac
{12}5NR-2\Omega B,\\
\gamma^2=&-\frac {6N(BP+A_{ 1.0})}{5A}+N_{ 1.0}+\frac 65NQ+ 2\Omega
A,
\endaligned
\end{equation}
in the case $B\ne 0$ it's components are given by the formulas
\begin{equation}\label{gamma2}
\aligned
\gamma^1=&-\frac {6N(AN-B_{ 0.1})}{5B}-N_{ 0.1} +\frac 65NR-2\Omega
B,\\
 \gamma^2=&-\frac {6AN(AS-B_{ 0.1})}{5B^2}+
\frac {18NAR}{5B}-
\frac {6N(A_{ 0.1}+B_{ 1.0})}{5B} +N_{ 1.0}-\frac
{12}5NQ+2\Omega A
\endaligned
\end{equation}

According to the Proposition all invariants (\ref{inv1}) are the constants.

\subsubsection{Case $I_1=const\ne 0$, $I_2=0$, $I_3=const$}

In the paper \cite{Kartak2} was described the Case $I_1=const$, $I_2=0$ from (\ref{inv1}).
Let us note that if $I_2=0$ (it means  $\Omega=0$ from (\ref{Omega1}), (\ref{Omega2})) then equation (\ref{eq}) by the change of variables (\ref{zam})
reduces into the form
\begin{equation}\label{model}
y''=P(x,y),
\end{equation}
proof is in the paper \cite{BordagBandle}. Let us remember that equations of the form (\ref{model})
were investigated in the works \cite{Lie}, \cite{Ovsyan}.

According to the paper \cite{Kartak2}, 
these equations  belong to one of four Types, depending on the value of the invariant  $I_1$, see Table. Here $t(x)$, $s(x)$ are the arbitrary functions
\begin{table}[ht]
{\begin{tabular}{|c|c|c|} 
\hline
Type & Equation &  $I_1$ \\
\hline
I & $y''=e^y+t(x)y+s(x)$ & $\frac 35$ \\
II & $y''=-\ln y+t(x)y+s(x)$  & $-\frac 9{10}$ \\
III & $y''=y(\ln y-1)+t(x)y+s(x)$ & $-\frac{12}{5}$ \\
IV & $y''=\frac{y^{C+2}}{(C+1)(C+2)}+t(x)y+s(x)$ &  $ \frac {3(C+5)}{5C}$, $C\ne 0,-1, -2,-5$\\
\hline
\end{tabular}}
\end{table}

Every Type splits into Subtypes.
Only subtypes I.1 and IV.1 satisfy the condition
$I_3=const$ from (\ref{inv1}). The Subtype I.1 is described by the following Theorem 2.

{\bf Theorem 2. }
 {\it These propositions are equivalent:
 \begin{enumerate}
 \item The equation (\ref{eq}) by the point transformations (\ref{zam}) reduces into the form $y''=e^y$;
 \item Equation (\ref{eq}) is from the Case 1.4 of Intermediate Degeneration, Subtype I.1: $F=0$ from (\ref{F}), $A\ne 0$ or $B\ne 0$ from (\ref{alpha}), $N\ne 0$ from (\ref{N}), $M\ne 0$ from (\ref{M1}), (\ref{M2}), $\Omega=0$ from (\ref{Omega1}), (\ref{Omega2}), $I_1=3/5$, $I_3=1/15$ from (\ref{inv1}).
 \end{enumerate}
 This equation has the 2-dimensional symmetry algebra.
 }
 
{\bf Example.} The Painlev\'e III equation with three zero parameters of the four  $a,\, b,\, c,\, d$
$$
y^{\prime\prime}=\frac     1y(y^{\prime})^2-\frac     1x
y^{\prime}+\frac 1x(ay^2+b)+cy^3+\frac dy
$$
 satisfies the conditions of Theorem 2, see \cite{Kartak3}.

{\bf Example.} The equations No. 6.76 (special case), 6.83 (special case), 6.110, 6.111 in the Handbook by E.Kamke \cite{Kamke} satisfy the conditions of Theorem 2
$$
\aligned
6.76\;\quad & y''=-\frac{a}{x}y'-be^y,\quad a=1,\, b\ne 0,\\
6.83\quad & y''=-\frac{a(e^y-1)}{x^2},\quad a=-2,\\
6.110 \quad & y''=\frac{{y'}^2}{y}-\frac 1y,\\
6.111 \quad & y''=\frac{{y'}^2}{y}+\frac 1y.\\
\endaligned
$$

The Subtype IV.1 is described by the  Theorem 3.

{\bf Theorem 3. }
 {\it These propositions are equivalent:
 \begin{enumerate}
 \item The equation (\ref{eq}) by the point transformations (\ref{zam})  reduces into the form  $y''=y^{c+2}/((c+1)(c+2))$, $c=const$, $c\ne -5,\,-2,\,-1,\,0$;
 \item The equation(\ref{eq}) is from the Case 1.4 of Intermediate Degeneration, Subtype IV.1: $F=0$ from (\ref{F}), $A\ne 0$ or $B\ne 0$ from (\ref{alpha}), $N\ne 0$ from (\ref{N}), $M\ne 0$ from (\ref{M1}), (\ref{M2}), $\Omega=0$ from (\ref{Omega1}), (\ref{Omega2}), $I_1=3(c+5)/(5c)$, $I_3=c(c+5)/(15(c+1)(c+2))$, $c=const$, $c\ne -5,\,-2,\,-1,\,0$ from (\ref{inv1}).
 \end{enumerate}
 This equation has the 2-dimensional symmetries algebra.
 }
 
{\bf Example.} The equations No. 6.7, 6.126 in the Handbook by E.Kamke \cite{Kamke} satisfy the conditions of Theorem 3
$$
\aligned
6.7 \quad & y''=ay^3,\; a\ne 0,\qquad I_1=\frac{18}{5},\quad  I_3=\frac 1{15},\quad c=1,\\
6.126 \quad & y y''+a(y'^2+1)=0,\, a\ne -3,\, -\frac 12,\, 0,\, 1,\\
&\phantom{.........}  I_1=-\frac{6(1+2a)}{5(3+a)},\; I_3=\frac {(3+a)(1+2a)}{15(a-1)},\; c=-\frac{a+3}{a+1}.\\
\endaligned
$$

\subsubsection{Case $I_1=const\ne 0$, $I_2=const\ne 0$, $I_3=const$}

As if $A$ and $B$ from (\ref{alpha}) are the components of the pseudovectorial field ${\mathbf \alpha}=(B,\,-A)$, then exists a special coordinate system, such that in the new coordinates the following conditions will be true $B=0$, $A=1$ (according to the rectification theorem). All possible change of variables (\ref{zam}) preserving the conditions $B=0$, $A=1$ has a transfer matrix of the form
$$
S=
\left(\begin{array}{cc}
	 \partial x /\partial \tilde  x& \partial x /\partial  \tilde y\\ \partial y /\partial  \tilde x & \partial  y /\partial \tilde y
	 \end{array}\right)=
\left(
\begin{array}{ll}
 s_{11} & 0 \\ 
s_{21}  &  1/{s_{11}^2}
 \end{array}
\right)
$$
Suppose that the equation (\ref{eq}) is already written in the such coordinate system.
Then   the condition $F=0$ from (\ref{F}) is true, so $S(x,y)=0.$ Hence
$$
N=RA\ne 0,\quad \Omega=Q_{1.0}-2R_{0.1}\ne 0,\quad I_1=-\frac{12}{5}-\frac{R_{1.0}}{R^2}=const\ne 0.
$$
Let us solve the equation $-R_{1.0}/{R^2}=c_1=const$, then
$R(x,y)={1}/({c_1y+f(x)}).$ 
By change of variables $x=\tilde x$, $y=\tilde y-f(\tilde x)/c_1$ we preserve the conditions $A=1$, $B=0$ and the coefficient $R$ will be the following $R(x,y)=1/(c_1y).$

Let us substitude $R$ into the formula for $B$ from (\ref{alpha}), we get
$$
Q_{0.2}+\frac{3Q_{0.1}}{c_1y}=0,\qquad\Rightarrow\qquad Q_{0.1}=g(x)y^{-\frac 3{c_1}}.
$$
Let us calculate the invariant $I_2$ from (\ref{inv1}):
$$
I_2=\frac{c_1yQ_{0.1}^2}{A}=c_2=const\ne 0, \qquad\Rightarrow\qquad c_1=6,\quad g(x)^2=\frac{2c_2}3.
$$
Hence $Q(x,y)=\pm\sqrt{6c_2y}/3.$
Let us calculate the invariant  $I_3$ from (\ref{inv1}), we get
$$
I_3=\frac{-72cy^2+P}{10y^2}=c_3=const,\qquad\Rightarrow\qquad P(x,y)=2(5c_3+36c_2)y^2.
$$
By the condition $A=1$ from (\ref{alpha}) follows that
$15c_3+110c_2=1.$ Then $c_3=1/15-22c_2/3.$

{\bf Theorem 4. }
 {\it These propositions are equivalent:
 \begin{enumerate}
 \item The equation (\ref{eq}) by the point transformations (\ref{zam}) reduces into the form  $y''=2(1-2k)y^2/3\pm\sqrt{6yk}y'+y'^2/(2y)$;
 \item The equation(\ref{eq}) is from the Case 1.2 of Intermediate Degeneration: $F=0$ from (\ref{F}), $A\ne 0$ or $B\ne 0$ from (\ref{alpha}), $N\ne 0$ from (\ref{N}), $M\ne 0$ from (\ref{M1}), (\ref{M2}), $\Omega\ne 0$ from(\ref{Omega1}), (\ref{Omega2}), $I_1=18/5$, $I_2=k=const\ne 0$, $I_3=1/15-22k/3$ from (\ref{inv1}).
 \end{enumerate}
 This equation has the 2-dimensional symmetries algebra.

 }
 
{\bf Example.} The equations No. 6.30, 6.174 in the Handbook by E.Kamke \cite{Kamke} satisfy the conditions of Theorem 4
$$
\aligned
6.30 \quad & y''+yy'-y^3=0,\qquad I_1=\frac{18}5,\quad I_2=\frac 1{20}=k,\quad I_3=-\frac 3{10},\\
6.174 \quad & xy y''-2xy'^2+(y+1)y'=0,\qquad I_1=\frac{18}5,\quad I_2=\frac 1{2}=k,\quad I_3=-\frac {18}{5}.\\
\endaligned
$$

\subsection{Seventh Case of the Intermediate Degeneration}

Seventh Case of the Intermediate Degeneration is characterized by the conditions $N=0$ from (\ref{N}),  $\Omega=0$ from (\ref{Omega1}), (\ref{Omega2}). All equations (\ref{eq}) from these case reduce into the form (\ref{model}), see \cite{Sharipov3}. 
In such form 
$$
A=P_{0.2}\ne 0,\quad B=0,\quad N=-\frac{A_{0.1}}{3}=0,\quad \Omega=0.
$$
Hence $A=f(x)$.
According to the paper  \cite{BordagBandle}, the most general  point transformations (\ref{zam})
preserving form (\ref{model}) are the following
\begin{equation}\label{form}
x=\alpha\int p^2(\tilde x) d(\tilde x)+\beta,\quad y=p(\tilde x)\tilde y+h(\tilde x).
\end{equation}
Here $\alpha$, $\beta$ are the constants, $p(\tilde x)$, $h(\tilde x)$ -- arbitrary functions.

So the direct and the inverse transfer matrices $S$ and $T$ are
$$
\aligned
S=&
\left(\begin{array}{cc}
	 \partial x /\partial \tilde  x& \partial x /\partial  \tilde y\\ \partial y /\partial  \tilde x & \partial  y /\partial \tilde y
	 \end{array}\right)=
\left(
\begin{array}{rr}
 \alpha p^2(\tilde x) & 0 \\ 
p_{1.0}(\tilde x)\tilde y+h_{0.1}(\tilde x)  & p(\tilde x)
 \end{array}
\right),\quad \det S=\alpha p^3(\tilde x),\\
 T=&S^{-1}=
\left(
\begin{array}{rr}
  \dfrac{1}{\alpha p^2(\tilde x)} & 0 \\ 
\dfrac{-p_{1.0}(\tilde x)\tilde y-h_{0.1}(\tilde x) }{\alpha p^3(\tilde x)} & \dfrac 1{ p(\tilde x)}
 \end{array}
\right),\qquad \qquad\quad\;\;\det T= \dfrac{1}{\alpha p^3(\tilde x)}.
\endaligned
$$
Under the transformations (\ref{form}) the pseudovectorial field ${\alpha} $ of the weight 2 (details are in the paper \cite{Sharipov3}) transforms as
$$
\left(\begin{array}{r}
\tilde B \\ -\tilde A
\end{array}\right)= \frac{T}{(\det T)^2}\left(\begin{array}{r}
 0 \\ -A
\end{array}\right)=
\left(\begin{array}{r}
0 \\ -\alpha p^5(\tilde x) A
\end{array}\right).
$$
So the components $A$ and $B$ transform by the rules
$
\tilde A(\tilde x,\tilde y)=\alpha^2p^5(\tilde x) A(x(\tilde x), y(\tilde x,\tilde y)),$ $ \tilde B=0.$
Choosing the appropriate function $p(\tilde x)$  we get
$\tilde A=\alpha^2p^5(\tilde x)f(x(\tilde x))=1$. 

So in the new variables the coefficient $P(x,y)=y^2/2+s(x)y+t(x). $
Let us calculate the pseudoinvariant $L$. The general formula for $L$ is from the paper \cite{Sharipov3}
\begin{equation}\label{L}
\aligned
L=&
 \theta^1\theta^2(\theta^1_{
1.0}-\theta^2_{ 0.1})+ (\theta^2)^2\theta^1_{
0.1}-(\theta^1)^2\theta^2_{ 1.0}-\\
& -P(\theta^1)^3-3Q(\theta^1)^2\theta^2-3R\theta^1(\theta^2)^2-S(\theta^2)^3-\frac 12\Theta^2,\quad\text{where}
\endaligned
\end{equation}
$$
\theta^1=\Theta_{ 0.1}-2\varphi_2\Theta, \qquad \theta^2=-\Theta_{
1.0}+2\varphi_1\Theta,\quad\text{and}\quad \Theta=\frac {\omega_1}A,\; A\ne 0, \qquad \Theta=\frac {\omega_2}B,\, B\ne 0,
$$
$\varphi_1$ and $\varphi_2$ are from the formula (\ref{phi}).
The components $\omega$  in the case $A\ne 0$ are equal to
$$
\aligned \omega_1=&\frac {12PR}{5A}-\frac {54}{25}\frac
{Q^2}A-\frac {P_{ 0.1}}A+ \frac {6Q_{ 1.0}}{5A}- \frac {PA_{
0.1}+BP_{ 1.0}+A_{ 2.0}}{5A^2}-\frac {2B_{ 1.0}P}{5A^2}+\\ &+
\frac {3QA_{ 1.0}-12PBQ}{25A^2}+
\frac {6B^2P^2+12BPA_{ 1.0}+6A_{ 1.0}^2}{25A^3},\\
\omega_2=& \frac {-5BP_{ 0.1}+6BQ_{ 0.1}+12RBP}{5A^2} -\frac
{54}{25}\frac {BQ^2}{A^2}+
\frac {6BA_{ 1.0}^2+6B^3P^2+12B^2A_{ 1.0}P}{25A^4}-\\ 
& -\frac {12B^2PQ}{25A^3}+ \frac {3BQA_{
1.0}}{25A^3}-\frac {2BB_{ 1.0}P+BA_{ 0.1}P+B^2P_{
1.0}+ BA_{ 2.0}}{5A^3}
\endaligned
$$
and in the case $B\ne 0$ write as
$$
\aligned \omega_1=& \frac {5AS_{ 1.0}-6AR_{ 0.1}+12QAS}{5B^2}
-\frac {54}{25}\frac {AR^2}{B^2}+
\frac {6AB_{ 0.1}^2+6A^3S^2-12A^2B_{ 0.1}S}{25B^4}-\\ 
&-\frac
{12A^2SR}{25B^3}+\frac {3ARB_{ 0.1}}{25B^3}+\frac {2AA_{ 0.1}S+AB_{
1.0}S+A^2S_{ 0.1}- AB_{ 0.2}}{5B^3},\\
\omega_2=&\frac {12SQ}{5B}-\frac {54}{25}\frac {R^2}B+\frac {S_{
1.0}}B- \frac {6R_{ 0.1}}{5B}+ \frac {SB_{ 1.0}+AS_{ 0.1}-B_{
0.2}}{5B^2}+\frac {2A_{ 0.1}S}{5B^2}-\\ &- \frac {3RB_{
0.1}+12SAR}{25B^2}+ \frac {6A^2S^2-12B_{ 0.1}AS+6B_{
0.1}^2}{25B^3}.
\endaligned
$$

For the equation $y''=y^2/2+s(x)y+t(x)$ the pseudoinvariant $L$ is given by the formula
$$
L=s''(x)-\frac{s^2(x)}{2}+t(x)=0,\quad\text{hence}\quad t(x)=\frac{s^2(x)}{2}-s''(x).
$$
By the change of variables
$x=\tilde x,$ $ y=\tilde y-s(\tilde x)$
this equation becomes 
$$
\tilde y''=\frac{{\tilde y}^2}2.
$$

{\bf Theorem 5. }
 {\it These propositions are equivalent:
 \begin{enumerate}
 \item The equation (\ref{eq}) by the point transformations (\ref{zam}) reduces into the form $$y''=y^2/2.$$
 \item The equation(\ref{eq}) is from the Seventh Case of Intermediate Degeneration: $F=0$ from (\ref{F}), $A\ne 0$ or $B\ne 0$ from(\ref{alpha}), $N= 0$ from (\ref{N}),  $\Omega=0$ from (\ref{Omega1}), (\ref{Omega2}), $L=0$ from (\ref{L}).
  \end{enumerate}
 This equation has the 2-dimensional  symmetry algebra.
 }

{\bf Example.} The equation No. 6.2 in the Handbook by E.Kamke \cite{Kamke} satisfies the conditions of Theorem 6
$$
\aligned
6.2 \quad & y''=6y^2.
\endaligned
$$

\subsection{General case}

The General Case is characterized by the condition $F\ne 0$ from (\ref{F}). In this case the basic invariants are calculated by the new rules, see \cite{Sharipov3} (the formula to the invariant $I_6$ is corrected here over the paper \cite{Sharipov3})
\begin{equation}
\begin{aligned}\label{inv2}
I_3&=\frac{B(HG_{1.0}-GH_{1.0})}{3F^9}-\frac{A(HG_{0.1}-GH_{0.1})}{3F^9}+
\frac{HF_{0.1}+GF_{1.0}}{3F^5}+\\
&+\frac{BG^2P}{3F^9}-\frac{(AG^2-2HBG)Q}{3F^9}+\frac{(BH^2-2HAG)R}{3F^9}-
\frac{AH^2S}{3F^9},\\
I_6&=\frac{H(AB_{0.1}-BA_{0.1})}{3F^7}+\frac{G(AB_{1.0}-BA_{1.0})}{3F^7}-
\frac{(AF_{0.1}-BF_{1.0})}{3F^3}-\\
&-\frac{GB^2P}{3F^7}-\frac{(HB^2-2GBA)Q}{3F^7}-\frac{(GA^2-2HBA)R}{3F^7}-
\frac{HA^2S}{3F^7},\\
I_7&=\frac{GHG_{1.0}-G^2H_{1.0}+H^2G_{0.1}-HGH_{0.1}+G^3P+3G^2HQ+3GH^2R+H^3S}{3F^{11}},\\
I_8&=\frac{G(AG_{1.0}+BH_{1.0})}{3F^9}+\frac{H(AG_{0.1}+BH_{0.1})}{3F^9}-
\frac{10(HF_{0.1}+GF_{1.0})}{3F^5}-\\
&-\frac{BG^2P}{3F^9}+\frac{(AG^2-2HBG)Q}{3F^9}-\frac{(BH^2-2HAG)R}{3F^9}+
\frac{AH^2S}{3F^9}.
\end{aligned}
\end{equation}

According to the Proposition, all invariants must be the constants.

Let for the equation (\ref{eq}) the conditions $B=0$, and $A=1$ from (\ref{alpha}) are true.
The most general type of point transformation that preserve these conditions is
\begin{equation}\label{z1}
x=p(\tilde x),\quad y=\frac{\tilde y}{{p'(\tilde x)}^2}+h(\tilde x)
\end{equation}
As if in this case $F^5=A^2S$, then $S\ne 0$. By the change of variables (\ref{z1})
the coefficient $S(x,y)$ is transformed by the rule $\tilde S=S(p'(\tilde x))^{-5}.$
 
\subsubsection{$S(x,y)=S(x)$} Let us the coefficient $S$ depends only on $x$, 
then choosing the appropriate function $p(\tilde x)$ we can get $S=1$, also conditions $B=0$, $A=1$ are saved.
Let us calculate the other invariants (\ref{inv2})
 $$
 \begin{aligned}
 I_3&= -3R_{0.1}-3Q+3R^2=c_3=const,\qquad \Rightarrow\qquad Q=-\frac {c_3}3-R_{0.1}+R^2;\\
 I_6&=0;\\
 I_7&= 9R_{1.0}+18RR_{0.1}+9P+9Rc_3-9R^3=c_7=const,\quad \Rightarrow\\
 &\phantom{...........................///////////////....}\quad P=\frac {c_7}9-R_{1.0}-2RR_{0.1}-Rc_3+R^3;\\
 I_8&= -c_3-3R_{0.1}=c_8=const\qquad \Rightarrow\qquad R=-(c_3+c_8)y/3+f(x).
 \end{aligned}
 $$
 Where $f(x)$ is an arbitrary function.
  Let us calculate the functions $Q$, $P$ and substitude them into the conditions $B$ and $A$ from (\ref{alpha}).
 $$
 \begin{aligned}
 B&=2c_8(c_3+c_8)/3=0,\\
 A&=(c_3+c_8)(-2yc_3c_8+6f(x)c_8+c_7-2yc_8^2)/9=1.
 \end{aligned}
 $$
So we get $c_8=0$, $c_7=9/c_3$.
 
Let us make a new change of variables (\ref{z1}): $x=\tilde x$, $y=\tilde y+a(\tilde x)$. Here the function
$a(\tilde x)$ satisfies the equation
$$
a'(\tilde x)=c_3a(\tilde x)/3-f(\tilde x).
$$
In the new variables the coeffitions will be
$$
\tilde R=-c_3\tilde y/3,\qquad \tilde Q=c_3^2{\tilde y}^2/9,\qquad \tilde P=1/c_3+c_3^2\tilde y/9-
c_3^3{\tilde y}^3/27.
$$

{\bf Theorem 6. }
 {\it These propositions are equivalent:
 \begin{enumerate}
 \item The equation (\ref{eq}) by the point transformations (\ref{zam}) reduces into the form $y''=y'^3-ky y'^2+k^2y^2y'/3+1/k+k^2y/9-k^3y^3/27$, $k=const$, $k\ne 0$;
 \item The equation(\ref{eq}) is from the Generl Case: $F\ne 0$ from (\ref{F}), with the conditions
 $I_3=k=const\ne 0$, $I_6=0$, $I_7=9/k$, $I_8=0$ from (\ref{inv2}).
  \end{enumerate}
 This equation has the 2-dimensional symmetry algebra.
 }

 {\bf Example.} The equation No. 6.130 in the Handbook by E.Kamke \cite{Kamke} in the case $c=b^2/(3(3+a))$, but $ab\ne 0$, $a\ne -3$ satisfies the conditions of Theorem 6 
 $$
 y''=-cy^3-byy'-\frac{ay^{\prime 2}}{y},\quad\text {where}\quad  b^2+c^2\ne 0.
 $$
The basic invariants from (\ref{inv2}) are
 $$
I_6=I_8=0,\qquad I_3=\sqrt[5]{\frac{3^6}{a^2(a+3)}}=k,\qquad I_7=\sqrt[5]{3^4a^2(a+3)}=\frac{9}{k}.
$$
This equation has the 2-dimensional algebra with the operators
 $$
 X_1=\frac{\partial}{\partial x},\qquad X_2=-x\frac{\partial}{\partial x}+y\frac{\partial}{\partial y}.
 $$

\subsubsection{$S(x,y)$ depends  on $x$ and $y$.} 

Let the coefficient $S$ depends  on $x$ and $y$ and the conditions $B=0$ and $A=1$ from (\ref{alpha}) are true. Let us calculate the invariant $I_6$ from (\ref{inv2})
$$
I_6=-\frac{S_{0.1}}{15S^{7/5}}=\frac{c_6}6=const\ne 0,\qquad \Rightarrow\qquad S=(c_6y+a(x))^{-5/2},$$
where $a(x)$ is an arbitrary function. By the following change of variables (\ref{z1}) 
$x=\tilde x$, $y=\tilde y-a(\tilde x)/c_6$
we get this function vanishing. Then
\begin{equation}\label{s}
S=\frac 1{(c_6y)^{5/2}}.
\end{equation}

As if $I_3=c_3=const$ from (\ref{inv2}) we can find
$$
Q=\frac 13\sqrt{c_6}\sqrt{y}(3c_6^2y^2(R^2-R_{0.1})-7c_6^2yR-c_3).
$$
From the formula $I_8=c_8=const$ from (\ref{inv2}) we get
$$
-\frac 92c_6^2yR-3c_6^2y^2R_{0.1}-c_3=c_8, \qquad \Rightarrow\qquad 
R=\frac{b(x)}{y^{3/2}}-\frac{2(c_3+c_8)}{3yc_6^2}.
$$
Let us make a point transformation  (\ref{z1}), preserving the function $S=(c_6y)^{-5/2}$
$$
y=\frac{\tilde y}{r^{\prime 2}(\tilde x)},\quad x=r(\tilde x),\quad\text{where}\quad b(r(\tilde x)) r^{\prime 2}(\tilde x)c_6^{5/2}-2r''(\tilde x)=c_6^{5/2}r'(\tilde x)
$$
After  the point transformation $R$ will be (we should not write the tildes on the variables)
\begin{equation}\label{r}
R=\frac{1}{y^{3/2}}-\frac{2(c_3+c_8)}{3yc_6^2},
\end{equation}
\begin{equation}\label{q}
Q=\frac{\sqrt{y}(c_6^2(5c_3+8c_8)+4(c_3+c_8)^2))}{9c_6^{3/2}}-
\frac{\sqrt{c_6}(5c_6^2+8(c_3+c_8))}{6}+\frac{c_6^{5/2}}{\sqrt{y}}
\end{equation}
From the formula $I_7=c_7=const$ from (\ref{inv2}) we get
\begin{equation}\label{p}
\begin{aligned}
P=&-\frac{8(c_3+c_8)^3y^2}{27c_6}-\frac{(4c_3^2+14c_3c_8+10c_8^2-c_7c_6)c_6y^2}{9}+\frac{4c_6(c_3+c_8)^2y^{3/2}}3+\\
&+\frac{c_6^3(5c_3+8c_8)y^{3/2}}{3}-\frac{c_6^3(3c_6^2+4(c_3+c_8))y}{2}+
c_6^5\sqrt{y}
\end{aligned}
\end{equation}

Now let us check the condition $B=0$, so
\begin{equation}\label{c7}
c_7=\frac{4c_8(c_3+c_8)}{c_6}-\frac{c_6(5c_3+8c_8)}{6}.
\end{equation}

Let us check the condition $A=1$. From this formula we can find $c_3$ as a solution of a quadratic equation
\begin{equation}\label{c3}
27c_6+(c_6^2-3c_8)(c_6^2(5c_3+8c_8)+4(c_3+c_8)^2)=0,
\end{equation}
let us note that $c_6^2-3c_8\ne 0$.

Let us rewrite formulas (\ref{s}), (\ref{r}), (\ref{q}), (\ref{p}) in the more appropriative form
$$
\begin{aligned}
S=&\frac{1}{ny^{5/2}},\quad R=\frac{1}{y^{3/2}}+\frac{a}{y},\quad 
Q=b\sqrt{y}+c +\frac{n}{\sqrt{y}},\\
P=&dy^2+ky^{3/2}+my+n^2\sqrt{y},\quad\text{where}
\end{aligned}
$$

\begin{equation}\label{co}
\begin{aligned}
n& ={c_6^{5/2}},\quad a=-\frac{2(c_3+c_8)}{3c_6^2},
\quad b= -\frac{3}{\sqrt{c_6}(c_6^2-3c_8)},\\  
c& =\frac{n(12a-5)}{6}, \quad d =\frac{bn(6a-1)}{6},\quad
k =3bn,\quad m=\frac{3n^2(2a-1)}2.
\end{aligned}
\end{equation}
Here $n$, $a$, $b$, $c$, $d$, $k$, $m$ are the constants. As if the parameters $c_3$, $c_6$ and $c_8$ satisfy to the condition (\ref{c3}), we  get a  relation between $n$, $a$ and $b$.
\begin{equation}\label{con}
18b(b-a^2n)+15nab-2bn-6=0.
\end{equation}

Let us make a point transformation $ y={\tilde y}^2,$ $ x=\tilde x$ 
 then the model equation will be the following (we should not write the tildes on the variables)
 \begin{equation}\label{t8}
 y''=\frac{4y'^3}{ny^3}+\frac{((6a-1)y+6)y'^2}{y^2}+\frac{3(by^2+cy+n)y'}{y}+
 \frac{1}2(dy^3+ky^2+my+n^2). 
   \end{equation}
 
 {\bf Theorem 7. }
 {\it These propositions are equivalent:
 \begin{enumerate}
 \item The equation (\ref{eq}) by the point transformations (\ref{zam}) reduces to the form (\ref{t8}) with the parameters (\ref{co}) such that the relation (\ref{con}) is true;
 \item The equation(\ref{eq}) is from the General Case: $F\ne 0$ из (\ref{F}), 
 invariants from (\ref{inv2}) are equal to $I_6=c_6/6=const\ne 0$, $I_8=c_8=const\ne {c_6^2}/{3}$, $I_7=c_7=const$  satisfy to the relation (\ref{c7}), $I_3=c_3=const$  satisfy to the relation (\ref{c3}).
  \end{enumerate}
 This equation has the 2-dimensional symmetries algebra  with the operators
   $$
X_1=\frac{\partial}{\partial x},\qquad X_2=e^{nx/2}\frac{\partial}{\partial x}-\frac{ny}{2}e^{nx/2}\frac{\partial}{\partial y}.
 $$}

 {\bf Example.} The equation No. 6.109 in the Handbook by E.Kamke \cite{Kamke} satisfies the conditions of Theorem 7
 $$
 y''=\frac{y'}{y}-\frac{y^{\prime 2}}{y}.
 $$
 The basic invariants from (\ref{inv2}) are equal to
  $$
 I_3=-\frac{37 \sqrt[5]{54}}{18},\quad 
 I_6=\frac{\sqrt[5]{648}}{9},\quad
 I_7=-\frac{55 \sqrt[5]{144}}{18},\quad
I_8= \frac{25 \sqrt[5]{54}}{18}. 
 $$
 The additional parameters are
 $$
 n=\frac{16\sqrt{3}}{3},\quad a=\frac 16,\quad b=\frac{\sqrt{3}}{3},\quad
 c=-\frac{8\sqrt{3}}{3},\quad d=0,\quad k=16,\quad m=-\frac{256}{3}.
 $$
 The canonical form will be
 $$
 y''= \frac{\sqrt{3}y'^3}{4y^3}+\frac{6y'^2}{y^2}+\frac{\sqrt{3}(4+y)^2y'}{y}+
 \frac{8}3(3y-4)(y-4). 
 $$

\section*{Acknowledgments}

 The work is partially supported by the Government of Russian Federation through Resolution No. 220, Agreement No. 11.G34.31.0042 and partially supported by the Russian Education and Science Ministry, Agreement No. 14.B37.21.0358.

\end{document}